\documentclass[12pt]{article}
\usepackage{color,amsfonts,amssymb}
\usepackage{amsfonts,epsf,amsmath}
\usepackage{latexsym}
\usepackage{graphicx}
\usepackage{array,float}
\usepackage{tikz}

\newtheorem{theorem}{Theorem}[section]
\newtheorem{lemma}[theorem]{Lemma}
\newtheorem{proposition}[theorem]{Proposition}

\newtheorem{conjecture}[theorem]{Conjecture}

\newcommand{\proof}{\noindent{\bf Proof.\ }}
\newcommand{\qed}{\hfill $\square$ \bigskip}
\newcommand{\cp}{\,\square\,}
\newcommand{\gp}{{\rm gp}}
\newcommand{\gpe}{{\rm gp_{e}}}

\begin{document}
\title{Edge general position sets in Fibonacci and Lucas cubes}

\author{
Sandi Klav\v zar$^{1, 2, 3}$
\and
Elif Tan$^{4}$
}

\date{}

\maketitle
\begin{center}
$^1$ Faculty of Mathematics and Physics,  University of Ljubljana, Slovenia\\
{\tt sandi.klavzar@fmf.uni-lj.si}\\
\medskip

$^2$ Faculty of Natural Sciences and Mathematics,  University of Maribor, Slovenia\\
\medskip

$^3$ Institute of Mathematics, Physics and Mechanics, Ljubljana, Slovenia\\
\medskip

$^4$ Department of Mathematics, Ankara University, Ankara, Turkey  \\
{\tt etan@ankara.edu.tr}\\

\end{center}

\begin{abstract}
A set of edges $X\subseteq E(G)$ of a graph $G$ is an edge general position set if no three edges from $X$ lie on a common shortest path in $G$. The cardinality of a largest edge general position set of $G$ is the edge general position number of $G$. In this paper edge general position sets are investigated in partial cubes. In particular it is proved that  the union of two largest $\Theta$-classes of a Fibonacci cube or a Lucas cube is a maximal edge general position set. 
\end{abstract}



\section{Introduction}

A set of vertices $S\subseteq V(G)$ of a graph $G = (V(G), E(G))$ is a {\em general position set} if no three vertices from $S$ lie on a common shortest path of $G$. Similarly, a set of edges $X\subseteq E(G)$ of $G$ is an {\em edge general position set} if no three edges from $X$ lie on a common shortest path. The cardinality of a largest general position set (resp.\ edge general position set) of $G$ is the {\em general position number} (resp.\ {\em edge general position number}) and denoted by $\gp(G)$ (resp.\ $\gpe(G)$). 

General position sets in graphs have already received a lot of attention. They were introduced in~\cite{manuel-2018, ullas-2016} and we refer to~\cite{AnaChaChaKlaTho, klavzar-2021, patkos-2020, tian-2021a, yao-2022} for a selection of further developments. On the other hand, the edge version of this concept has been introduced only recently in~\cite{manuel-2022}. In this paper we continue this line of the research. 

To determine the general position number of hypercubes turns out to be a very difficult problem, cf.~\cite{korner-1995}. On the other hand, a closed formula for the edge general position number of hypercubes has been determined in~\cite{manuel-2022}. Combining the facts that the edge general position number is doable on hypercubes and that hypercubes form the cornerstone of the class of partial cubes, we focus in this paper on the edge general position number of two important and interesting families of partial cubes, Fibonacci cubes and Lucas cubes. The first of these two classes of graphs was introduced in~\cite{hsu-1993} as a model for interconnection networks. In due course, these graphs have found numerous applications elsewhere and are also extremely interesting in their own right. Lucas cubes, introduced in~\cite{munarini-2001}, form a class of graphs which naturally symmetrises the Fibonacci cubes and also have many interesting properties. The state of research up to 2013 on these classes of graphs (and some additional related ones) is summarised in the survey paper~\cite{klavzar-2013}, the following list of papers is a selection from subsequent research~\cite{egecioglu-2021a, egecioglu-2021b, gravier-2015, ilic-2017, mollard-2021, mollard-2022, savitha-2020, saygi-2018, saygi-2019, taranenko-2013, vesel-2015}. 
 
The rest of this paper is organized as follows. In the next section we define the concepts discussed in this paper, introduce the required notation, and recall a known result. In  Section~\ref{sec:partial-cubes} we discuss partial cubes and the interdependence of their edge general position sets and $\Theta$-classes. In Section~\ref{sec:Fibonacci-Lucas} we prove that the union of two largest $\Theta$-classes of a Fibonacci cube or a Lucas cube is always a maximal edge general position set. We conjecture that for Fibonacci cubes these sets are also maximum general position sets and show that this is not the case for Lucas cubes. 

\section{Preliminaries}
\label{sec:prelim}

Unless stated otherwise, graphs considered in this paper are connected. The path of order $n$ is denoted by $P_n$. The {\em Cartesian product} $G\cp H$ of graphs $G$ and $H$ has vertices $V(G)\times V(H)$ and edges $(g,h)(g',h')$, where either $g=g'$ and $hh'\in E(H)$, or $h=h'$ and $gg'\in E(G)$. The {\em $r$-dimensional hypercube} $Q_r$, $r\ge 1$, is a graph with $V(Q_r) = \{0, 1\}^r$, and there is an edge between two vertices if and only if they differ  in exactly one coordinate. That is, if $x = (x_1, \ldots, x_r)$ and $y = (y_1, \ldots, y_r)$ are vertices of $Q_r$, then $xy\in E(Q_r)$ if and only if there exists  $j\in [r] = \{1,\ldots, r\}$ such that $x_j \neq y_j$ and $x_i = y_i$ for every $i\ne j$. $Q_r$, $r\ge 2$, can also be described as the Cartesian product $Q_{r-1}\cp P_2$. 

The {\em distance} $d_G(u, v)$ between vertices $u$ and $v$ of a graph $G = (V(G), E(G))$ is the number of edges on a shortest $u, v$-path. A subgraph $H$ of a graph $G$ is {\em isometric} if $d_H(x, y) = d_G(x, y)$ holds for every pair of vertices $x, y$ of $H$. We also say that $H$ is {\em isometrically embedded} into $G$. Isometric subgraphs of hypercubes are known as {\em partial cubes}. 

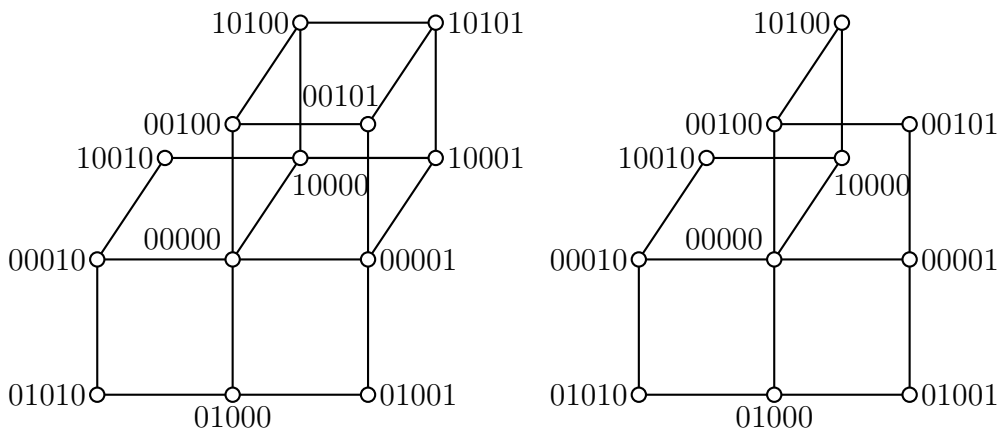
\begin{figure}[ht!]
\begin{center}
\begin{tikzpicture}[scale=0.9,style=thick]
\def\vr{3pt}

\begin{scope}[xshift=0cm, yshift=0cm]
\path (0,0) coordinate (01010);
\path (2,0) coordinate (01000);
\path (4,0) coordinate (01001);
\path (0,2) coordinate (00010);
\path (2,2) coordinate (00000);
\path (4,2) coordinate (00001);
\path (1,3.5) coordinate (10010);
\path (3,3.5) coordinate (10000);
\path (5,3.5) coordinate (10001);
\path (2,4) coordinate (00100);
\path (4,4) coordinate (00101);
\path (3,5.5) coordinate (10100);
\path (5,5.5) coordinate (10101);
\draw (01010) -- (01000) --(01001) -- (00001) -- (00000) -- (00010) -- (10010) -- (10000) -- (10001) -- (10101) --(10100) -- (10000) -- (10001);
\draw (01010) -- (00010); 
\draw (01000) -- (00000); 
\draw (00000) -- (10000); 
\draw (00001) -- (10001); 
\draw (00100) -- (10100); 
\draw (00001) -- (00101); 
\draw (00000) -- (00100); 
\draw (00101) -- (10101); 
\draw (00100) -- (00101); 
\draw (00000)  [fill=white] circle (\vr);
\draw (10000)  [fill=white] circle (\vr);
\draw (01000)  [fill=white] circle (\vr);
\draw (00100)  [fill=white] circle (\vr);
\draw (00010)  [fill=white] circle (\vr);
\draw (00001)  [fill=white] circle (\vr);
\draw (10100)  [fill=white] circle (\vr);
\draw (10010)  [fill=white] circle (\vr);
\draw (10001)  [fill=white] circle (\vr);
\draw (10101)  [fill=white] circle (\vr);
\draw (01010)  [fill=white] circle (\vr);
\draw (01001)  [fill=white] circle (\vr);
\draw (00101)  [fill=white] circle (\vr);
\draw[left] (00000)++(0,0.3) node {$00000$}; 
\draw[right] (10000)++(-0.3,-0.4) node {$10000$}; 
\draw[below] (01000)++(0,0.0) node {$01000$}; 
\draw[left] (00100)++(0,0.0) node {$00100$}; 
\draw[left] (00010)++(0,0.0) node {$00010$}; 
\draw[right] (00001)++(0,0.0) node {$00001$}; 
\draw[left] (10100)++(0,0.0) node {$10100$}; 
\draw[left] (10010)++(0,0.0) node {$10010$}; 
\draw[right] (10001)++(0,0.0) node {$10001$}; 
\draw[left] (01010)++(0,0.0) node {$01010$}; 
\draw[right] (01001)++(0,0.0) node {$01001$}; 
\draw[above] (00101)++(-0.4,0.1) node {$00101$}; 
\draw[right] (10101)++(0,0.0) node {$10101$}; 
\end{scope}

\begin{scope}[xshift=8cm, yshift=0cm]
\path (0,0) coordinate (01010);
\path (2,0) coordinate (01000);
\path (4,0) coordinate (01001);
\path (0,2) coordinate (00010);
\path (2,2) coordinate (00000);
\path (4,2) coordinate (00001);
\path (1,3.5) coordinate (10010);
\path (3,3.5) coordinate (10000);
\path (5,3.5) coordinate (10001);
\path (2,4) coordinate (00100);
\path (4,4) coordinate (00101);
\path (3,5.5) coordinate (10100);
\path (5,5.5) coordinate (10101);
\draw (01010) -- (01000) --(01001) -- (00001) -- (00000) -- (00010) -- (10010) -- (10000);
\draw  (10100) -- (10000);
\draw (01010) -- (00010); 
\draw (01000) -- (00000); 
\draw (00000) -- (10000); 
\draw (00100) -- (10100); 
\draw (00001) -- (00101); 
\draw (00000) -- (00100); 
\draw (00100) -- (00101); 
\draw (00000)  [fill=white] circle (\vr);
\draw (10000)  [fill=white] circle (\vr);
\draw (01000)  [fill=white] circle (\vr);
\draw (00100)  [fill=white] circle (\vr);
\draw (00010)  [fill=white] circle (\vr);
\draw (00001)  [fill=white] circle (\vr);
\draw (10100)  [fill=white] circle (\vr);
\draw (10010)  [fill=white] circle (\vr);
\draw (01010)  [fill=white] circle (\vr);
\draw (01001)  [fill=white] circle (\vr);
\draw (00101)  [fill=white] circle (\vr);
\draw[left] (00000)++(0,0.3) node {$00000$}; 
\draw[right] (10000)++(-0.3,-0.4) node {$10000$}; 
\draw[below] (01000)++(0,0.0) node {$01000$}; 
\draw[left] (00100)++(0,0.0) node {$00100$}; 
\draw[left] (00010)++(0,0.0) node {$00010$}; 
\draw[right] (00001)++(0,0.0) node {$00001$}; 
\draw[left] (10100)++(0,0.0) node {$10100$}; 
\draw[left] (10010)++(0,0.0) node {$10010$}; 
\draw[left] (01010)++(0,0.0) node {$01010$}; 
\draw[right] (01001)++(0,0.0) node {$01001$}; 
\draw[right] (00101)++(0,0) node {$00101$}; 
\end{scope}

\end{tikzpicture}
\end{center}
\caption{$\Gamma_5$ and $\Lambda_5$}
\label{fig:Gamma_5-Lambda_5}
\end{figure}

A {\em Fibonacci string} of length $n\ge 1$ is a binary string that contains no consecutive 1s. The {\em Fibonacci cube} $\Gamma_n$, $n\ge 1$, is the graph whose vertices are all Fibonacci strings of length $n$, two vertices being adjacent if they differ in a single coordinate. $\Gamma_n$ can be equivalently defined as an induced subgraph of $Q_n$ obtained from $Q_n$ by removing all the vertices that contain at least one pair of consecutive 1s. Further, the {\em Lucas cube} $\Lambda_n$, $n\ge 1$, is obtained from $\Gamma_n$ by removing the vertices that start and end with $1$. See Fig.~\ref{fig:Gamma_5-Lambda_5} for $\Gamma_5$ and $\Lambda_5$ and note that the latter is obtained from the former by removing the vertices $10001$ and $10101$. 

It is well-known that the order of $\Gamma_{n}$ is $F_{n+2}$, where $F_n$ are the {\em Fibonacci numbers} defined by the recurrence $F_{n+2} = F_{n+1} + F_{n}$, $n\ge 0$, with the initial terms $F_0 = 0$ and $F_1 = 1$. Also, the order of $\Lambda_{n}$ is $L_n$, where $L_n$ are the {\em Lucas numbers} defined by the same recurrence relation with the initial terms $L_0 = 2$ and $L_1 = 1$. 

To complete the preliminaries, we recall the following inequality on Fibonacci numbers.  

\begin{lemma} {\rm \cite[Corollary]{attanasov-2003}}  
\label{lem:Hosoya}
If $n$ is a positive integer and $0 \le i\le \lfloor n/2\rfloor$, then $F_{n}\geq F_{i}F_{n-i+1}$.
\end{lemma}

\section{On edge general position sets in partial cubes}
\label{sec:partial-cubes}

In this section we recall several results on the edge general position sets in partial cubes
from~\cite{manuel-2022} and derive some new results. This will motivate us to consider edge general position sets in Fibonacci cubes and in Lucas cubes in the next section.

Let $G$ be a graph. Then we say that edges $xy$ and $uv$ of $G$ are in the {\em Djokovi\'c-Winkler relation} $\Theta$ if $d_G(x,u) + d_G(y,v) \not= d_G(x,v) + d_G(y,u)$~\cite{djok-73,wink-84}. A connected graph $G$ is a partial cube if and only if $G$ is bipartite and $\Theta$ is transitive~\cite{wink-84}. It follows that $\Theta$ partitions the edge set of a partial cube into $\Theta$-{\em classes}. Moreover, if $G$ is a partial cube isometrically embedded into $Q_n$ such that for each $i\in [n]$ there exits an edge $xy\in E(G)\subseteq E(Q_n)$ with $x_i\ne y_i$, then $G$ contains exactly $n$ $\Theta$-classes. We will denote them by $\Theta_1(G), \ldots, \Theta_n(G)$, where $\Theta_i(G)$, $i\in [n]$, consists of the edges of $G$ which differ in coordinate $i$.   

We first recall the following result. 

\begin{lemma} {\rm \cite[Lemma~3.1]{manuel-2022}}
\label{lem:two-classes}
If $G$ is a partial cube embedded into $Q_n$, then $\Theta_i(G) \cup \Theta_j(G)$ is an edge general position set of $G$. 
\end{lemma}

Using Lemma~\ref{lem:two-classes} it was proved in~\cite{manuel-2022} that $\gpe(Q_r) = 2^r$. It was also proved that $\gpe(P_r\cp P_r) = 4r - 8$ for $r\ge 4$, see~\cite[Theorem 4.1]{manuel-2022}. Note however, that $|\Theta_i(P_r\cp P_r)| = r$ for each $i\in [2r-2]$. Hence Lemma~\ref{lem:two-classes} only yields $\gpe(P_r\cp P_r) \ge 2r$ which is arbitrary away from the optimal value. Moreover, if $r\ge 5$, then by~\cite[Theorem~4.2]{manuel-2022} a largest edge general position set of $P_r\cp P_r$ is unique and is not a union of some $\Theta$-classes. On the other hand, we can have large edge general position sets which are the union of many $\Theta$-classes as the following result asserts, where by an {\em end block} of a graph $G$ we mean a block of $G$ which contains exactly one cut vertex. 

To prove the next proposition, we recall the following auxiliary result. 

\begin{lemma} {\rm \cite[Lemma~6.4]{klavzar-2013b}}
\label{lem-subgraph} 
Let $H$ be an isometric subgraph of $G$ and let $e$ and $f$ be edges from different blocks of $H$. Then $e$ is not in relation $\Theta$ with $f$ in $G$.
\end{lemma}

\begin{proposition}
\label{prop:block}
Let $B_1,\ldots, B_k$ be the end blocks of a partial cube $G$ and for $i\in [k]$ let $\Theta^{i}(G)$ be an arbitrary $\Theta$-class of $G$ with an edge in $B_i$. Then $\bigcup_{i\in [k]} \Theta^{i}(G)$ is an edge general position set of $G$. 
\end{proposition}

\proof
By Lemma~\ref{lem-subgraph}, $\Theta^{i}(G) \subseteq E(B_i)$ and thus $\Theta^{i}(G)$ also forms a $\Theta$-class of $B_i$. Consider now an arbitrary shortest path $P$ of $G$ and suppose it contains an edge $e_i$ of some $\Theta^{i}(G)$. Then by the above, $e_i$ is the only edge of $P$ from $\Theta^{i}(G)$. If $P$ contains an edge $e_j$ from some other $\Theta^{j}(G)$, then also $e_j$ is the only edge of $P$ from $\Theta^{j}(G)$. Moreover, in this case, $E(P) \cap \bigcup_{i\in [k]} \Theta^{i}(G) = \{e_i, e_j\}$ because all the edges of $P$ which do not lie in $B_i \cup B_j$ are from blocks which are not end blocks. 
\qed

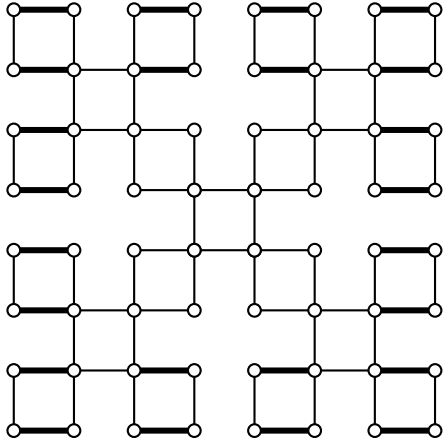
\begin{figure}[ht!]
\begin{center}
\begin{tikzpicture}[scale=0.8,style=thick]
\def\vr{3pt}

\draw (-2,2) -- (-2,3) -- (-1,3) -- (-1,2) -- (-2,2); 
\draw (2,2) -- (2,3) -- (3,3) -- (3,2) -- (2,2); 
\draw (-2,-1) -- (-2,-2) -- (-1,-2) -- (-1,-1) -- (-2,-1); 
\draw (2,-1) -- (2,-2) -- (3,-2) -- (3,-1) -- (2,-1); 

\begin{scope}[xshift=0cm, yshift=0cm]
\path (0,0) coordinate (a);
\path (1,0) coordinate (b);
\path (1,1) coordinate (c);
\path (0,1) coordinate (d);
\draw (a) -- (b) -- (c) -- (d) -- (a); 
\draw (a)  [fill=white] circle (\vr);
\draw (b)  [fill=white] circle (\vr);
\draw (c)  [fill=white] circle (\vr);
\draw (d)  [fill=white] circle (\vr);
\end{scope}

\begin{scope}[xshift=-3cm, yshift=1cm]
\path (0,0) coordinate (a);
\path (1,0) coordinate (b);
\path (1,1) coordinate (c);
\path (0,1) coordinate (d);
\draw (a) -- (b) -- (c) -- (d) -- (a); 
\draw[line width=0.8mm] (a) -- (b); 
\draw[line width=0.8mm] (c) -- (d); 
\draw (a)  [fill=white] circle (\vr);
\draw (b)  [fill=white] circle (\vr);
\draw (c)  [fill=white] circle (\vr);
\draw (d)  [fill=white] circle (\vr);
\end{scope}

\begin{scope}[xshift=-1cm, yshift=1cm]
\path (0,0) coordinate (a);
\path (1,0) coordinate (b);
\path (1,1) coordinate (c);
\path (0,1) coordinate (d);
\draw (a) -- (b) -- (c) -- (d) -- (a); 
\draw (a)  [fill=white] circle (\vr);
\draw (b)  [fill=white] circle (\vr);
\draw (c)  [fill=white] circle (\vr);
\draw (d)  [fill=white] circle (\vr);
\end{scope}

\begin{scope}[xshift=1cm, yshift=1cm]
\path (0,0) coordinate (a);
\path (1,0) coordinate (b);
\path (1,1) coordinate (c);
\path (0,1) coordinate (d);
\draw (a) -- (b) -- (c) -- (d) -- (a); 
\draw (a)  [fill=white] circle (\vr);
\draw (b)  [fill=white] circle (\vr);
\draw (c)  [fill=white] circle (\vr);
\draw (d)  [fill=white] circle (\vr);
\end{scope}

\begin{scope}[xshift=3cm, yshift=1cm]
\path (0,0) coordinate (a);
\path (1,0) coordinate (b);
\path (1,1) coordinate (c);
\path (0,1) coordinate (d);
\draw (a) -- (b) -- (c) -- (d) -- (a); 
\draw[line width=0.8mm] (a) -- (b); 
\draw[line width=0.8mm] (c) -- (d); 
\draw (a)  [fill=white] circle (\vr);
\draw (b)  [fill=white] circle (\vr);
\draw (c)  [fill=white] circle (\vr);
\draw (d)  [fill=white] circle (\vr);
\end{scope}

\begin{scope}[xshift=-3cm, yshift=3cm]
\path (0,0) coordinate (a);
\path (1,0) coordinate (b);
\path (1,1) coordinate (c);
\path (0,1) coordinate (d);
\draw (a) -- (b) -- (c) -- (d) -- (a); 
\draw[line width=0.8mm] (a) -- (b); 
\draw[line width=0.8mm] (c) -- (d); 
\draw (a)  [fill=white] circle (\vr);
\draw (b)  [fill=white] circle (\vr);
\draw (c)  [fill=white] circle (\vr);
\draw (d)  [fill=white] circle (\vr);
\end{scope}

\begin{scope}[xshift=-1cm, yshift=3cm]
\path (0,0) coordinate (a);
\path (1,0) coordinate (b);
\path (1,1) coordinate (c);
\path (0,1) coordinate (d);
\draw (a) -- (b) -- (c) -- (d) -- (a); 
\draw[line width=0.8mm] (a) -- (b); 
\draw[line width=0.8mm] (c) -- (d); 
\draw (a)  [fill=white] circle (\vr);
\draw (b)  [fill=white] circle (\vr);
\draw (c)  [fill=white] circle (\vr);
\draw (d)  [fill=white] circle (\vr);
\end{scope}

\begin{scope}[xshift=1cm, yshift=3cm]
\path (0,0) coordinate (a);
\path (1,0) coordinate (b);
\path (1,1) coordinate (c);
\path (0,1) coordinate (d);
\draw (a) -- (b) -- (c) -- (d) -- (a); 
\draw[line width=0.8mm] (a) -- (b); 
\draw[line width=0.8mm] (c) -- (d); 
\draw (a)  [fill=white] circle (\vr);
\draw (b)  [fill=white] circle (\vr);
\draw (c)  [fill=white] circle (\vr);
\draw (d)  [fill=white] circle (\vr);
\end{scope}

\begin{scope}[xshift=3cm, yshift=3cm]
\path (0,0) coordinate (a);
\path (1,0) coordinate (b);
\path (1,1) coordinate (c);
\path (0,1) coordinate (d);
\draw (a) -- (b) -- (c) -- (d) -- (a); 
\draw[line width=0.8mm] (a) -- (b); 
\draw[line width=0.8mm] (c) -- (d); 
\draw (a)  [fill=white] circle (\vr);
\draw (b)  [fill=white] circle (\vr);
\draw (c)  [fill=white] circle (\vr);
\draw (d)  [fill=white] circle (\vr);
\end{scope}

\begin{scope}[xshift=-3cm, yshift=-1cm]
\path (0,0) coordinate (a);
\path (1,0) coordinate (b);
\path (1,1) coordinate (c);
\path (0,1) coordinate (d);
\draw (a) -- (b) -- (c) -- (d) -- (a); 
\draw[line width=0.8mm] (a) -- (b); 
\draw[line width=0.8mm] (c) -- (d); 
\draw (a)  [fill=white] circle (\vr);
\draw (b)  [fill=white] circle (\vr);
\draw (c)  [fill=white] circle (\vr);
\draw (d)  [fill=white] circle (\vr);
\end{scope}

\begin{scope}[xshift=-1cm, yshift=-1cm]
\path (0,0) coordinate (a);
\path (1,0) coordinate (b);
\path (1,1) coordinate (c);
\path (0,1) coordinate (d);
\draw (a) -- (b) -- (c) -- (d) -- (a); 
\draw (a)  [fill=white] circle (\vr);
\draw (b)  [fill=white] circle (\vr);
\draw (c)  [fill=white] circle (\vr);
\draw (d)  [fill=white] circle (\vr);
\end{scope}

\begin{scope}[xshift=1cm, yshift=-1cm]
\path (0,0) coordinate (a);
\path (1,0) coordinate (b);
\path (1,1) coordinate (c);
\path (0,1) coordinate (d);
\draw (a) -- (b) -- (c) -- (d) -- (a); 
\draw (a)  [fill=white] circle (\vr);
\draw (b)  [fill=white] circle (\vr);
\draw (c)  [fill=white] circle (\vr);
\draw (d)  [fill=white] circle (\vr);
\end{scope}

\begin{scope}[xshift=3cm, yshift=-1cm]
\path (0,0) coordinate (a);
\path (1,0) coordinate (b);
\path (1,1) coordinate (c);
\path (0,1) coordinate (d);
\draw (a) -- (b) -- (c) -- (d) -- (a); 
\draw[line width=0.8mm] (a) -- (b); 
\draw[line width=0.8mm] (c) -- (d); 
\draw (a)  [fill=white] circle (\vr);
\draw (b)  [fill=white] circle (\vr);
\draw (c)  [fill=white] circle (\vr);
\draw (d)  [fill=white] circle (\vr);
\end{scope}

\begin{scope}[xshift=-3cm, yshift=-3cm]
\path (0,0) coordinate (a);
\path (1,0) coordinate (b);
\path (1,1) coordinate (c);
\path (0,1) coordinate (d);
\draw (a) -- (b) -- (c) -- (d) -- (a); 
\draw[line width=0.8mm] (a) -- (b); 
\draw[line width=0.8mm] (c) -- (d); 
\draw (a)  [fill=white] circle (\vr);
\draw (b)  [fill=white] circle (\vr);
\draw (c)  [fill=white] circle (\vr);
\draw (d)  [fill=white] circle (\vr);
\end{scope}

\begin{scope}[xshift=-1cm, yshift=-3cm]
\path (0,0) coordinate (a);
\path (1,0) coordinate (b);
\path (1,1) coordinate (c);
\path (0,1) coordinate (d);
\draw (a) -- (b) -- (c) -- (d) -- (a); 
\draw[line width=0.8mm] (a) -- (b); 
\draw[line width=0.8mm] (c) -- (d); 
\draw (a)  [fill=white] circle (\vr);
\draw (b)  [fill=white] circle (\vr);
\draw (c)  [fill=white] circle (\vr);
\draw (d)  [fill=white] circle (\vr);
\end{scope}

\begin{scope}[xshift=1cm, yshift=-3cm]
\path (0,0) coordinate (a);
\path (1,0) coordinate (b);
\path (1,1) coordinate (c);
\path (0,1) coordinate (d);
\draw (a) -- (b) -- (c) -- (d) -- (a); 
\draw[line width=0.8mm] (a) -- (b); 
\draw[line width=0.8mm] (c) -- (d); 
\draw (a)  [fill=white] circle (\vr);
\draw (b)  [fill=white] circle (\vr);
\draw (c)  [fill=white] circle (\vr);
\draw (d)  [fill=white] circle (\vr);
\end{scope}

\begin{scope}[xshift=3cm, yshift=-3cm]
\path (0,0) coordinate (a);
\path (1,0) coordinate (b);
\path (1,1) coordinate (c);
\path (0,1) coordinate (d);
\draw (a) -- (b) -- (c) -- (d) -- (a); 
\draw[line width=0.8mm] (a) -- (b); 
\draw[line width=0.8mm] (c) -- (d); 
\draw (a)  [fill=white] circle (\vr);
\draw (b)  [fill=white] circle (\vr);
\draw (c)  [fill=white] circle (\vr);
\draw (d)  [fill=white] circle (\vr);
\end{scope}

\end{tikzpicture}
\end{center}
\caption{An edge general position set}
\label{fig:pattern}
\end{figure}

Note that Proposition~\ref{prop:block} implies that the set of leaves of a tree $T$ forms an edge general position set. For another example of an edge general position which is the union of many $\Theta$-classes see Fig.~\ref{fig:pattern}. 

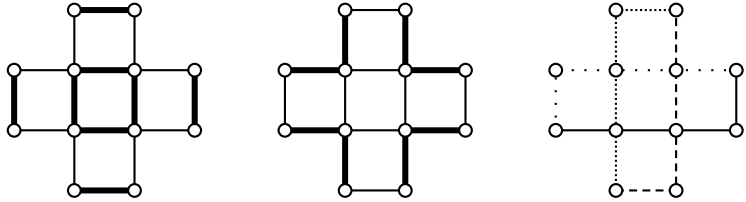
\begin{figure}[ht!]
\begin{center}
\begin{tikzpicture}[scale=0.8,style=thick]
\def\vr{3pt}

\begin{scope}[xshift=0cm, yshift=0cm]
\path (1,0) coordinate (x1);
\path (2,0) coordinate (x2);
\path (0,1) coordinate (x3);
\path (1,1) coordinate (x4);
\path (2,1) coordinate (x5);
\path (3,1) coordinate (x6);
\path (0,2) coordinate (x7);
\path (1,2) coordinate (x8);
\path (2,2) coordinate (x9);
\path (3,2) coordinate (x10);
\path (1,3) coordinate (x11);
\path (2,3) coordinate (x12);
\draw (x1) -- (x11) -- (x12) -- (x2) -- (x1); 
\draw (x3) -- (x6) -- (x10) -- (x7) -- (x3); 
\draw[line width=0.8mm] (x3) -- (x7); 
\draw[line width=0.8mm] (x4) -- (x8); 
\draw[line width=0.8mm] (x5) -- (x9); 
\draw[line width=0.8mm] (x6) -- (x10); 
\draw[line width=0.8mm] (x11) -- (x12); 
\draw[line width=0.8mm] (x8) -- (x9); 
\draw[line width=0.8mm] (x4) -- (x5); 
\draw[line width=0.8mm] (x1) -- (x2); 
\foreach \i in {1,...,12}
{
\draw (x\i)  [fill=white] circle (\vr);
}
\end{scope}

\begin{scope}[xshift=4.5cm, yshift=0cm]
\path (1,0) coordinate (x1);
\path (2,0) coordinate (x2);
\path (0,1) coordinate (x3);
\path (1,1) coordinate (x4);
\path (2,1) coordinate (x5);
\path (3,1) coordinate (x6);
\path (0,2) coordinate (x7);
\path (1,2) coordinate (x8);
\path (2,2) coordinate (x9);
\path (3,2) coordinate (x10);
\path (1,3) coordinate (x11);
\path (2,3) coordinate (x12);
\draw (x1) -- (x11) -- (x12) -- (x2) -- (x1); 
\draw (x3) -- (x6) -- (x10) -- (x7) -- (x3); 
\draw[line width=0.8mm] (x3) -- (x4); 
\draw[line width=0.8mm] (x7) -- (x8); 
\draw[line width=0.8mm] (x8) -- (x11); 
\draw[line width=0.8mm] (x9) -- (x12); 
\draw[line width=0.8mm] (x9) -- (x10); 
\draw[line width=0.8mm] (x5) -- (x6); 
\draw[line width=0.8mm] (x2) -- (x5); 
\draw[line width=0.8mm] (x1) -- (x4); 
\foreach \i in {1,...,12}
{
\draw (x\i)  [fill=white] circle (\vr);
}
\end{scope}

\begin{scope}[xshift=9cm, yshift=0cm]
\path (1,0) coordinate (x1);
\path (2,0) coordinate (x2);
\path (0,1) coordinate (x3);
\path (1,1) coordinate (x4);
\path (2,1) coordinate (x5);
\path (3,1) coordinate (x6);
\path (0,2) coordinate (x7);
\path (1,2) coordinate (x8);
\path (2,2) coordinate (x9);
\path (3,2) coordinate (x10);
\path (1,3) coordinate (x11);
\path (2,3) coordinate (x12);
\draw[densely dotted] (x1) -- (x11) -- (x12); 
\draw[densely dashed] (x1) -- (x2) -- (x12);
\draw[loosely dotted] (x3) -- (x7) -- (x10); 
\draw (x3) -- (x6) -- (x10);
\foreach \i in {1,...,12}
{
\draw (x\i)  [fill=white] circle (\vr);
}
\end{scope}

\end{tikzpicture}
\end{center}
\caption{Two largest edge general position sets of a partial cube and a cover of it by shortest paths}
\label{fig:example}
\end{figure}

To finish this section consider the partial cube $G$ from Fig.~\ref{fig:example}. In the left figure, the union of its two largest $\Theta$-classes forms an edge general position set of cardinality $8$. In the middle figure, the union of four $\Theta$-classes also forms an edge general position set of $G$ of cardinality $8$. Finally, since we can cover the edges of $G$ by four shortest paths as shown in the right figure, we have $\gpe(G)\le 8$ which means that both indicated sets are largest edge general position sets.  

\section{On edge general position sets in Fibonacci and Lucas cubes}
\label{sec:Fibonacci-Lucas}

Fibonacci cubes and Lucas cubes are partial cubes~\cite{klavzar-2005}. Thus all the results and comments of the previous section can be applied to them. The cardinality of the $\Theta$-classes of Fibonacci cubes was independently determined in~\cite{klavzar-2007, taranenko-2007}, and the cardinality of the $\Theta$-classes of Lucas cubes in~\cite{klavzar-2007}. These results read as follows.  

\begin{proposition}
\label{prop:size-of-Theta-classes}
(i) If $n\geq 1$ and $i\in [n]$, then $|\Theta_i(\Gamma_n)| = F_{i}F_{n-i+1}$.

(ii) If $n\geq 1$ and $i\in [n]$, then $|\Theta_i(\Lambda_n)| = F_{n-1}$.
\end{proposition}

To find large edge general position sets in $\Gamma_n$ we can apply Lemma~\ref{lem:two-classes}. For this sake we first answer the question for which $i$ and $j$ the value $|\Theta_i(\Gamma_n) \cup \Theta_j(\Gamma_n)|$ is maximum. 

\begin{proposition}
If $n\ge 2$, then 
$$2F_n = \max \{|\Theta_i(\Gamma_n)| + |\Theta_j(\Gamma_n)|:\ i, j\in [n], i\ne j\}.$$
\end{proposition}

\proof
Set $M = \max \{|\Theta_i(\Gamma_n)| + |\Theta_j(\Gamma_n)|:\ i, j\in [n], i\ne j\}$. Using Proposition~\ref{prop:size-of-Theta-classes} and Lemma~\ref{lem:Hosoya} we can then estimate as follows: 
\begin{align*}
M & = \max \{|\Theta_i(\Gamma_n)| + |\Theta_j(\Gamma_n)|:\ i, j\in [n], i\ne j\} \\
& = \max \{F_iF_{n-i+1} + F_jF_{n-j+1}:\ i, j\in [n], i\ne j\} \\
& \le \max \{F_n + F_n:\ i, j\in [n], i\ne j\} \\ 
& = 2F_n\,.
\end{align*}
On the other hand, $|\Theta_1(\Gamma_n)| + |\Theta_n(\Gamma_n)| = F_n + F_n$, hence we can conclude that $M = 2F_n$. 
\qed

\begin{theorem}
\label{thm:maximal-Gamma}
If $n\ge 2$, then $\Theta_1(\Gamma_n) \cup \Theta_n(\Gamma_n)$ is a maximal edge general position set of $\Gamma_n$. Moreover, $\gpe(\Gamma_n)\ge 2F_n$.   
\end{theorem}

\proof
To prove the first assertion, we will use the fact that a shortest path can contain at most one edge from $\Theta_1(\Gamma_n)$ and at most one edge from $\Theta_n(\Gamma_n)$, cf.~\cite[Lemma~11.1]{hammack-2011}. Hence we only need to prove that no edge can be added to $\Theta_1(\Gamma_n) \cup \Theta_n(\Gamma_n)$ in order to keep the edge general position property. 

The statement of the theorem clearly holds for $\Gamma_2$, and can be easily verified for $\Gamma_3$ and $\Gamma_4$. In the rest we may thus assume that $n\ge 5$. Consider an arbitrary edge $e = uv\in \Theta_i(\Gamma_n)$, where $2\le i\le n-1$. We may without loss of generality assume that $u_i = 0$ and $v_i = 1$. We need to show that $e$ lies on some shortest path that contains one edge from $\Theta_1(\Gamma_n)$ and one edge from $\Theta_n(\Gamma_n)$. We distinguish the following two cases. 

\medskip\noindent
{\bf Case 1}: $i \in \{2, n-1\}$. \\
Suppose first that $i = 2$. In this case $u = 000 \ldots$ and $v = 010\ldots$ If $u_n = v_n = 1$, then the following path 
\begin{align*}
y & = 1\,0\,0\, \ldots \, 0\,0 \\
x & = 0\,0\,0\, \ldots \, 0\,0 \\
u & = 0\,0\,0\, \ldots \, 0\,1 \\
v & = 0\,1\,0\, \ldots \, 0\,1 
\end{align*}
is a shortest $x,y$-path in $\Gamma_n$ that contains $yx\in \Theta_1(\Gamma_n)$ and $xu\in \Theta_n(\Gamma_n)$. If $u_n = v_n = 0$, then we consider two subcases. In the first one, $u = 000 \ldots 00$ and $v = 010\ldots 00$. Then the following shortest path 
\begin{align*}
x & = 1\,0\,0\, \ldots \, 0\,0 \\
u & = 0\,0\,0\, \ldots \, 0\,0 \\
v & = 0\,1\,0\, \ldots \, 0\,0 \\
y & = 0\,1\,0\, \ldots \, 0\,1 
\end{align*}
contains $xu\in \Theta_1(\Gamma_n)$ and $vy\in \Theta_n(\Gamma_n)$. In the second subcase we consider  $u = 000 \ldots 10$ and $v = 010\ldots 10$, in which case we have  $u = 000 \ldots 010$ and $v = 010\ldots 010$. Then the path 
\begin{align*}
x & = 1\,0\,0\, \ldots \, 0\,1\,0 \\
u & = 0\,0\,0\, \ldots \, 0\,1\,0 \\
v & = 0\,1\,0\, \ldots \, 0\,1\,0 \\
y & = 0\,1\,0\, \ldots \, 0\,0\,0 \\
z & = 0\,1\,0\, \ldots \, 0\,0\,1 
\end{align*}
is a shortest path in $\Gamma_n$ and contains $xu\in \Theta_1(\Gamma_n)$ and $yz\in \Theta_n(\Gamma_n)$. For instance, if $n=5$, then the path constructed is: $x=10010$, $u=00010$, $v=01010$, $y=01000$, $z=01001$. 

We have thus considered all the subcases when $i = 2$. By the symmetry of Fibonacci strings, the case $i = n-1$ can be done analogously. 

\medskip\noindent
{\bf Case 2}: $2 < i < n-1$. \\
In this case we have $u = \ldots 000 \ldots$ and $v = \ldots 010 \ldots$ Assume first that $u_1= v_1 = 1$ and $u_n = v_n = 1$, so that $u = 10 \ldots 010\ldots  01$ and $v = 10 \ldots 000 \ldots 01$. Then the path 
\begin{align*}
x & = 0\,0\, \ldots \, 0\,1\,0 \ldots 0\, 1\\
v & = 1\,0\, \ldots \, 0\,1\,0 \ldots 0\, 1\\
u & = 1\,0\, \ldots \, 0\,0\,0 \ldots 0\, 1\\
y & = 1\,0\, \ldots \, 0\,0\,0 \ldots 0\, 0
\end{align*}
is a shortest path in $\Gamma_n$ and contains edges $xv\in \Theta_1(\Gamma_n)$ and $uy\in \Theta_n(\Gamma_n)$. For instance, if $n=5$, then the path constructed is: $x=00101$, $v=10101$, $u=10001$, $y=10000$.

If $u$ and $v$ start and end by $00$, we simply change the first and the last bit to construct a required shortest path. Assume next that $u = 01 \ldots 000\ldots  10$ and $v = 01 \ldots 010 \ldots 10$. Because $v$ starts and ends with $0$, and contains at least three $1$s, in this case we have $n\ge 7$. Then the path 
\begin{align*}
x & = 1\,0\, \ldots \, 0\,1\,0 \ldots 1\, 0\\
x' & = 0\,0\, \ldots \, 0\,1\,0 \ldots  1\, 0\\
v & = 0\,1\, \ldots \, 0\,1\,0 \ldots 1\, 0\\
u & = 0\,1\, \ldots \, 0\,0\,0 \ldots  1\, 0\\
y & = 0\,1\, \ldots \, 0\,0\,0 \ldots  0\, 0\\
y' & = 0\,1\, \ldots \, 0\,0\,0 \ldots  0\, 1
\end{align*}
is a shortest path in $\Gamma_n$ which contains $xx'\in \Theta_1(\Gamma_n)$ and $yy'\in \Theta_n(\Gamma_n)$. The final cases when $u$ and $v$ begin by $0$ and end by $1$ (or the other way around) are done by combining the above paths. 

By Proposition~\ref{prop:size-of-Theta-classes}(i), $|\Theta_1(\Gamma_n)| = |\Theta_n(\Gamma_n)| = F_n$, hence $\gpe(\Gamma_n)\ge |\Theta_1(\Gamma_n)| + |\Theta_n(\Gamma_n)| = 2F_n$.   
\qed

With a lot of effort, we can further prove that for any $i$ and $j$, the set $\Theta_i(\Gamma_n)\cup \Theta_j(\Gamma_n)$ is a maximal edge general position set. However, since $\Theta_1(\Gamma_n)\cup \Theta_n(\Gamma_n)$ is a largest such a set, we omit the long case analysis here. 

Based on Theorem~\ref{thm:maximal-Gamma} we wonder whether $\Theta_1(\Gamma_n) \cup \Theta_n(\Gamma_n)$ is not only a maximal edge general position set of $\Gamma_n$ but also a maximum edge general position set. While we have no answer in general, we next show that this is true up to dimension $n\le 5$. For $n\le 4$ this can be easily checked, and for $n=5$ we have the following. 

\begin{proposition}
\label{prop:up-to-5}
$\gpe(\Gamma_5) = 10$. 
\end{proposition}

\begin{figure}[ht!]
\begin{center}
\begin{tikzpicture}[scale=1.0,style=thick]
\def\vr{3pt}

\begin{scope}[xshift=0cm, yshift=0cm]
\path (0,0) coordinate (01010);
\path (2,0) coordinate (01000);
\path (4,0) coordinate (01001);
\path (0,2) coordinate (00010);
\path (2,2) coordinate (00000);
\path (4,2) coordinate (00001);
\path (1,3.5) coordinate (10010);
\path (3,3.5) coordinate (10000);
\path (5,3.5) coordinate (10001);
\path (2,4) coordinate (00100);
\path (4,4) coordinate (00101);
\path (3,5.5) coordinate (10100);
\path (5,5.5) coordinate (10101);
\draw[densely dotted] (10101) -- (10100) -- (10000) -- (10010) -- (00010); 
\draw[loosely dashed] (10101) -- (10001) -- (00001) -- (01001) -- (01000) -- (01010); 
\draw[densely dashed] (10101) -- (00101) -- (00001) -- (00000) -- (00010) -- (01010); 
\draw (10001) -- (10000) -- (00000) -- (01000); 
\draw[line width=0.8mm] (00100) -- (00000); 
\draw[line width=0.8mm] (00100) -- (10100); 
\draw[line width=0.8mm] (00100) -- (00101); 
\draw (00000)  [fill=white] circle (\vr);
\draw (10000)  [fill=white] circle (\vr);
\draw (01000)  [fill=white] circle (\vr);
\draw (00100)  [fill=white] circle (\vr);
\draw (00010)  [fill=white] circle (\vr);
\draw (00001)  [fill=white] circle (\vr);
\draw (10100)  [fill=white] circle (\vr);
\draw (10010)  [fill=white] circle (\vr);
\draw (10001)  [fill=white] circle (\vr);
\draw (10101)  [fill=white] circle (\vr);
\draw (01010)  [fill=white] circle (\vr);
\draw (01001)  [fill=white] circle (\vr);
\draw (00101)  [fill=white] circle (\vr);
\draw (1.8,3) node {$e$}; 
\draw (2.3,4.8) node {$e'$}; 
\draw (3.5,4.25) node {$e''$}; 
\end{scope}

\end{tikzpicture}
\end{center}
\caption{Four shortest paths in $\Gamma_5$}
\label{fig:four-paths}
\end{figure}
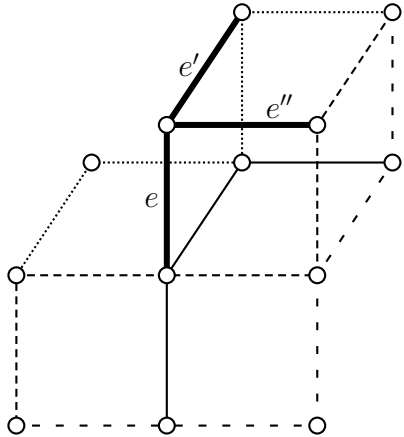

\proof
By Theorem~\ref{thm:maximal-Gamma} we have $\gpe(\Gamma_5) \ge 10$. Consider now the four paths in $\Gamma_5$ as indicated in Fig.~\ref{fig:four-paths} and note that each of them is a shortest path. The only edges not contained in one of these four paths are $e$, $e'$, and $e''$, see Fig.~\ref{fig:four-paths} again.  

Let $X$ be an arbitrary edge general position set of $\Gamma_5$. Since each of the four paths from Fig.~\ref{fig:four-paths} is a shortest path, each of then can contain at most two edges from $X$ and thus $|X| \le 4\cdot 2 + 3 = 11$. Supposing that $|X| = 11$, we must have $e, e', e''\in X$. However, by inspection we can now infer that if $e, e', e''\in X$, then there are only $5$ more additional edges that could possibly lie in $X$, hence if $e, e', e''\in X$, then actually $|X| \le 8$ holds. 
\qed

Based on Theorem~\ref{thm:maximal-Gamma} and Proposition~\ref{prop:up-to-5} we pose: 

\begin{conjecture}
\label{conj:Gamma_n}
If $n\ge 2$, then $\gpe(\Gamma_n) = 2F_n$. 
\end{conjecture}

For the Lucas cubes we have the following result parallel to Theorem~\ref{thm:maximal-Gamma}.

\begin{theorem}
\label{thm:maximal-Lambda}
If $n\ge 4$, then $\Theta_1(\Lambda_n) \cup \Theta_n(\Lambda_n)$ is a maximal edge general position set of $\Lambda_n$. Moreover, $\gpe(\Lambda_n)\ge 2F_{n-1}$.   
\end{theorem}

\proof
We proceed parallel with the proof of Theorem~\ref{thm:maximal-Gamma}. More precisely, we need to show that no edge can be added to $\Theta_1(\Lambda_n) \cup \Theta_n(\Lambda_n)$ in order to keep the edge general position property. Then all the paths constructed in the proof of Theorem~\ref{thm:maximal-Gamma} contain no vertex which would start and end with $1$, hence the same paths are suitable also for the present proof. The only exception appears to be the path as constructed in the first subcase of Case 2. However, this subcase is not relevant in the present proof because the vertices $u$ and $v$ are not vertices of $\Lambda_n$ and hence we need not consider them here. Finally, by Proposition~\ref{prop:size-of-Theta-classes}(ii), $|\Theta_1(\Lambda_n)| = |\Theta_n(\Lambda_n)| = F_{n-1}$, hence $\gpe(\Lambda_n)\ge |\Theta_1(\Lambda_n)| + |\Theta_n(\Lambda_n)| = 2F_{n-1}$.   
\qed

By Theorem~\ref{thm:maximal-Lambda}, $\gpe(\Lambda_5)\ge 6$. The following result then comes as a surprise. 

\begin{proposition}
\label{prop:Lambda_5}
$\gpe(\Lambda_5) = 7$. 
\end{proposition}

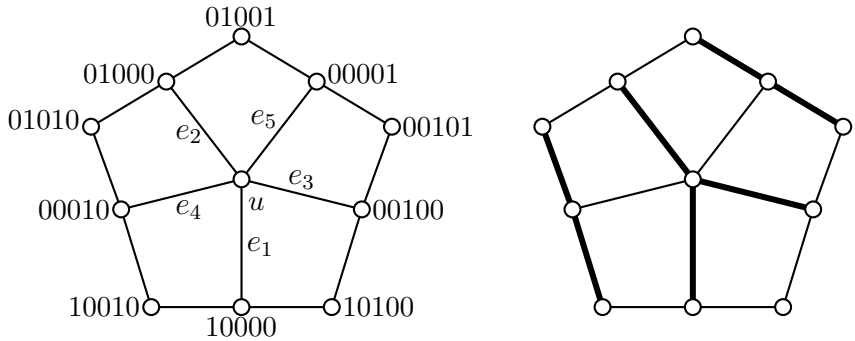
\begin{figure}[ht!]
\begin{center}
\begin{tikzpicture}[scale=1.0,style=thick]
\def\vr{3pt}

\begin{scope}[xshift=0cm, yshift=0cm]
\path (0,-0.3) coordinate (u);
\path (0,-2.0) coordinate (x1);
\path (1.6,-0.7) coordinate (x2);
\path (1,1) coordinate (x3);
\path (-1,1) coordinate (x4);
\path (-1.6,-0.7) coordinate (x5);
\path (1.2,-2.0) coordinate (y1);
\path (2.0,0.4) coordinate (y2);
\path (0,1.6) coordinate (y3);
\path (-2.0,0.4) coordinate (y4);
\path (-1.2,-2.0) coordinate (y5);
\foreach \i in {1,...,5}
{
\draw (u) -- (x\i); 
}
\draw (x1) -- (y1) -- (x2) -- (y2) -- (x3) -- (y3) -- (x4) --(y4) -- (x5) -- (y5) -- (x1); 

\draw (u)  [fill=white] circle (\vr);
\foreach \i in {1,...,5}
{
\draw (x\i)  [fill=white] circle (\vr);
\draw (y\i)  [fill=white] circle (\vr);
}
{\small 
\draw[below] (u)++(0.2,-0.1) node {$u$}; 
\draw[below] (x1)++(0.0,0.0) node {$10000$}; 
\draw[right] (x2)++(0.0,0.0) node {$00100$}; 
\draw [right](x3)++(0.0,0.1) node {$00001$}; 
\draw [left](x4)++(0.0,0.1) node {$01000$}; 
\draw[left] (x5)++(0.0,0.0) node {$00010$}; 
}
\draw (0.25,-1.2) node {$e_1$};
\draw (0.8,-0.3) node {$e_3$};
\draw (0.3,0.5) node {$e_5$};
\draw (-0.7,0.3) node {$e_2$};
\draw (-0.7,-0.7) node {$e_4$};
{\small 
\draw[right] (y1)++(0.0,0.0) node {$10100$}; 
\draw[right] (y2)++(0.0,0.0) node {$00101$}; 
\draw [above](y3)++(0.0,0.0) node {$01001$}; 
\draw [left](y4)++(0.0,0.1) node {$01010$}; 
\draw[left] (y5)++(0.0,0.0) node {$10010$}; 
}
\end{scope}

\begin{scope}[xshift=6cm, yshift=0cm]
\path (0,-0.3) coordinate (u);
\path (0,-2.0) coordinate (x1);
\path (1.6,-0.7) coordinate (x2);
\path (1,1) coordinate (x3);
\path (-1,1) coordinate (x4);
\path (-1.6,-0.7) coordinate (x5);
\path (1.2,-2.0) coordinate (y1);
\path (2.0,0.4) coordinate (y2);
\path (0,1.6) coordinate (y3);
\path (-2.0,0.4) coordinate (y4);
\path (-1.2,-2.0) coordinate (y5);
\foreach \i in {1,...,5}
{
\draw (u) -- (x\i); 
}
\draw (x1) -- (y1) -- (x2) -- (y2) -- (x3) -- (y3) -- (x4) --(y4) -- (x5) -- (y5) -- (x1); 
\draw[line width=0.8mm] (u) -- (x1); 
\draw[line width=0.8mm] (u) -- (x2); 
\draw[line width=0.8mm] (u) -- (x4); 
\draw[line width=0.8mm] (y5) -- (x5); 
\draw[line width=0.8mm] (x5) -- (y4); 
\draw[line width=0.8mm] (x3) -- (y3); 
\draw[line width=0.8mm] (x3) -- (y2); 
\draw (u)  [fill=white] circle (\vr);
\foreach \i in {1,...,5}
{
\draw (x\i)  [fill=white] circle (\vr);
\draw (y\i)  [fill=white] circle (\vr);
}
\end{scope}

\end{tikzpicture}
\end{center}
\caption{$\Lambda_5$ and its largest edge general position set}
\label{fig:Lambda_5-again}
\end{figure}

\proof
Let $X$ be an arbitrary edge general position set of $\Lambda _{5}$ and let $u=00000$. Let $F = \left\{ e_{1},e_{2},e_{3},e_{4},e_{5}\right\}$ be the set of the edges incident to $u$,  where $e_{1}=\{u,10000\}$, $e_{2} = \{u,01000\}$, $e_{3} = \{u, 00100\}$, $e_{4} = \{u,  00010\}$, and $e_{5} = \{u, 00001\}$. See Fig.~\ref{fig:Lambda_5-again}. 

Since $F$ is an edge general position set, we consider the following cases.

\medskip\noindent
{\bf Case 1}: $\left\vert X\cap F\right\vert =5$. \\
In this case no edge from $E(\Lambda _{5})\backslash F$ can be added to $X$, so $|X| = 5$.

\medskip\noindent
\textbf{Case 2}: $\left\vert X\cap F\right\vert =4$. \\
We may without loss of generality assume that $F=\{e_1, e_3, e_5, e_2\}$. Then only the edges $\{01010, 00010\}$ and $\{10010, 00010\}$ can be added to $X$, hence in this case we have $|X|\le 6$. 

\medskip\noindent
\textbf{Case 3}: $\left\vert X\cap F\right\vert =3$.\\
By the symmetry of $\Lambda_5$ it suffices to distinguish the following two subcases.

(i) $X$ contains three consecutive edges, say $X\cap F=\left\{ e_{1},e_{3},e_{5}\right\}$. Then at most two edges of $E(\Lambda _{5})\backslash F$ can be added to $X$.

(ii) $X$ does not contain three consecutive edges, say $X\cap F=\left\{ e_{1},e_{2},e_{3}\right\}$. Then four edges of $E(\Lambda _{5})\backslash F$ can be added to $X$ as shown in the right-hand side of Fig.~\ref{fig:Lambda_5-again}. So in this case $|X| = 7$. 

\medskip\noindent
\textbf{Case 4}: $\left\vert X\cap F\right\vert =2$. \\
In this case, no matter whether the edges from $X\cap F$ are consecutive or not, we can easily check that at most four edges of $E(\Lambda _{5})\backslash F$ can be added to $X$.

\medskip\noindent
\textbf{Case 5}: $\left\vert X\cap F\right\vert =1$. \\
We may assume that $F = \{e_1\}$. Then considering the three subcases based on whether $\{10000, 10010\}$ and $\{10000, 10100\}$ lie in $X$, we get that in every case $|X| \le 6$. 

\medskip\noindent
\textbf{Case 6:}  $\left\vert X\cap F\right\vert =0$. \\
In this case we observe that at most every second edge from the outer $10$-cycle can lie in $X$, hence $|X|\le 5$. 
\qed

Note that the proof of Proposition~\ref{prop:Lambda_5} also implies that a largest edge general position set of $\Lambda_5$ is unique up to symmetry.

\section*{Acknowledgements}

This work has been supported by T\"{U}B\.{I}TAK and the Slovenian Research Agency under grant numbers 122N184 and BI-TR/22-24-20, respectively. Sandi Klav\v{z}ar also acknowledges the financial support from the Slovenian Research Agency (research core funding P1-0297 and projects J1-2452 and N1-0285).

\section*{Declaration of interests}

The authors declare that they have no known competing financial interests or personal relationships that could have appeared to influence the work reported in this paper.

\section*{Data availability}

Our manuscript has no associated data.

\end{document}